\documentclass[12pt,twoside]{article}

\usepackage{amsmath,amssymb,amsthm,amscd,graphicx,color,accents,mathrsfs}

\usepackage[pagebackref=false,colorlinks=true,linkcolor=blue,citecolor=blue,urlcolor=blue]{hyperref}
\usepackage{url}

\numberwithin{equation}{section}

\theoremstyle{plain}
\newtheorem{theorem}{Theorem}[section]
\newtheorem{proposition}[theorem]{Proposition}
\newtheorem{lemma}[theorem]{Lemma}

\newtheorem{remark}[theorem]{Remark}

\newtheorem{example}[theorem]{Example}

\makeatletter
\newcommand*\bigcdot{\mathpalette\bigcdot@{.8}}
\newcommand*\bigcdot@[2]{\mathbin{\vcenter{\hbox{\scalebox{#2}{{\hskip 1pt}$\m@th#1\bullet$}}}}}
\makeatother

\newcommand{\Pro}{\noindent\textit{Proof.}\ \ }

\DeclareMathOperator\E {\mathbb{E}}
\def\G {{{\hskip 1pt}\mathscr{G}}}
\DeclareMathOperator\N {\mathbb{N}}
\def\P {\mathcal{P}}
\def\S {\mathcal{S}}
\DeclareMathOperator\R  {\mathbb{R}}
\DeclareMathOperator\C  {\mathbb{C}}
\def\Pr {\mathbb{P}}
\DeclareMathOperator\tr {\mathop{\rm{tr}\hskip 1pt}}
\DeclareMathOperator\Dom {\mathop{\rm{Dom}\hskip 0.2pt}}

\newcommand{\eqdist}{\overset{\mathcal{L}}{=}}
\newcommand{\gedist}{\overset{\mathcal{L}}{\ge}}
\newcommand{\ledist}{\overset{\mathcal{L}}{\le}}

\def \dd {\,\textrm{d}}

\def\diag   {{\rm{diag}}}

\begin{document}

\title{\Large\textbf{Hoffmann-J{\o}rgensen Inequalities for Random Walks on the Cone of Positive Definite Matrices}}

\author{
{Armine Bagyan}\thanks{Department of Statistics, Pennsylvania State University, University Park, PA 16802, U.S.A. \ E-mail address: aub171@psu.edu}
\ {and Donald Richards}\thanks{Department of Statistics, Pennsylvania State University, University Park, PA 16802, U.S.A. \ E-mail address: richards@stat.psu.edu
\endgraf
\ $^\dag$Corresponding author.
}}

\date{\today}

\maketitle

\begin{abstract}
We consider random walks on the cone of $m \times m$ positive definite matrices, where the underlying random matrices have orthogonally invariant distributions on the cone and the Riemannian metric is the measure of distance on the cone.  By applying results of Khare and Rajaratnam ({\it Ann.~Probab.}, 45 (2017), 4101--4111), we obtain inequalities of Hoffmann-J{\o}rgensen type for such random walks on the cone.  In the case of the Wishart distribution $W_m(a,I_m)$, with index parameter $a$ and matrix parameter $I_m$, the identity matrix, we derive explicit and computable bounds for each term appearing in the Hoffmann-J{\o}rgensen inequalities.  



\smallskip
\noindent
{\em Running head}: Hoffmann-J{\o}rgensen Inequalities for Positive Definite Matrices.
\end{abstract}

\section{Introduction}
\label{sec:intro}

In this paper, we consider inequalities of Hoffmann-J{\o}rgensen type for random walks on the cone of positive definite matrices.  These inequalities are obtained by adapting results of Khare and Rajaratnam \cite{Khare} obtained in the broader setting of metric semigroups.  In studying these inequalities on matrix cones, we are motivated by the appearance of random samples of positive definite matrices in numerous fields, including:~statistical inference on Riemannian manifolds \cite{Healy,Kim,Kim_Koo}, microwave engineering \cite[p.~156 ff.]{Terras}, diffusion tensor imaging \cite{Dryden,Jian}, financial time series \cite{Asai}, wireless communication systems \cite{Siriteanu}, polarimetric radar imaging \cite{Anfinsen}, factor analysis \cite{Browne}, and goodness-of-fit testing \cite{Hadjicosta2}.  

The resulting generalized Hoffmann-J{\o}rgensen inequalities on $\P_m$, the cone of $m\times m$ positive definite (symmetric) matrices, are derived under the assumption that the underlying random matrices have an orthogonally invariant distribution on the cone.  Specializing to the case of the Wishart distribution $W_m(a,I_m)$, with index parameter $a$ and identity matrix parameter $I_m$, we derive explicit and computable bounds for each term appearing in the Hoffmann-J{\o}rgensen inequalities.  Although we restrict our attention in this paper to $\P_m$, our methods apply more generally to the cone of Hermitian positive definite matrices and to abstract symmetric cones \cite{Ding,Faraut}.  

In Section \ref{sec_matrixcase}, we provide notation and preliminaries for random walks on $\P_m$ and for the natural action of $GL(m,\R)$, the group of nonsingular $m \times m$ matrices, on the random walks.  We provide in this section some necessary properties of the Riemannian metric, which is perhaps the most prominent metric on $\P_m$, and we show in Proposition \ref{prop:invariancemetrics} that the probability distribution of the distance between orthogonally invariant random matrices is preserved under orthogonally invariant random walks. 

In Section \ref{sec_gen_HJ_ineqs} we derive from a result of Khare and Rajaratnam \cite[Theorem A]{Khare} a generalization, to the cone $\P_m$, of the classical Hoffmann-J{\o}rgensen inequalities.  This result, which is stated in Theorem \ref{thm_HJ1}, provides bounds on the probability distribution of the distance between $X_n$, the $n$th ``step'' in an invariant random walk from the initial matrix $I_m$, the $m \times m$ identity matrix.  

In order to apply an Hoffmann-J{\o}rgensen inequality in a practical setting, it is necessary to calculate or obtain an upper bound for each term appearing in the inequality; this will be seen to be a non-trivial problem for the random walks considered here.  In Section \ref{sec_wishart}, we consider the special case in which the random walk is generated by random samples from the Wishart distribution $W_m(a,I_m)$, with index parameter $a$ and matrix parameter $I_m$, the identity matrix; we are particularly interested in Wishart-distributed random walks because of their appearance in the literature on diffusion tensor imaging, factor analysis, and financial volatility.  Under the Wishart assumption, we obtain explicit, computable bounds for each term in the Hoffmann-J{\o}rgensen inequalities in Theorem \ref{thm_HJ1}.  

Finally, in Appendix \ref{appendix}, we establish the submartingale properties of some random variables appearing in the generalized Hoffmann-J{\o}rgensen inequalities and show that the corresponding Kolmogorov inequalities reduce to Markov's inequalities for those variables.

\section{Random matrices and random walks on \texorpdfstring{$\boldsymbol{\P_m}$}{P(m,R)}}
\label{sec_matrixcase}

We denote by $G$ the general linear group $GL(m,\R)$ of all real, $m \times m$, nonsingular matrices.  The group $G$ {\it acts} on $\P_m$ {\it via} the correspondence, $g \cdot A = gAg'$, where $g \in G$, $A \in \P_m$, and $g'$ denotes the transpose of $g$.  It is elementary to verify that this correspondence is a {\it group action}, i.e., $g_1g_2 \cdot A = g_1 \cdot (g_2 \cdot A)$, for all $g_1, g_2 \in G$.

The group action also is transitive:~given $A, B \in \P_m$, there exists $g \in G$ such that $g \cdot A = B$; simply choose $g = B^{1/2} A^{-1/2}$.  Alternatively, since $A, B \in \P_m$, there exist $g_1, g_2 \in G$ such that $A = g_1g_1'$, $B = g_2 g_2'$; then $g = g_2 g_1^{-1}$ satisfies $g \cdot A = B$.  Denote by $I_m$ the identity matrix in $G$.  Under this group action, the isotropy group of $I_m$ is $O(m)$, the maximal compact subgroup of $m \times m$ orthogonal matrices.  Thus, the homogeneous space $G/O(m)$ can be identified with $\P_m$ through the correspondence $g\cdot k  = (gk)(gk)' = gg'$ for all $k \in O(m)$.

Whenever two random entities $X$ and $Y$ have the same probability distribution, we will write $X \eqdist Y$.  Then a random matrix $X \in \P_m$ is said to be {\it orthogonally invariant} if $X \eqdist kXk'$ for any $k \in O(m)$.  Examples of orthogonally invariant distributions are the well-known Wishart distribution $W_m(\alpha,I_m)$ \cite[Chapter 7]{Anderson} and the multivariate beta distribution \cite[p.~377]{Anderson}.  Similarly, a function $f :\P_m \to \R$ is said to be \textit{orthogonally invariant} if $f(X) = f(kXk')$ for all $k \in O(m)$ and $X \in \P_m$.  

If the random matrix $X = gg' \in \P_m$ is orthogonally invariant, where $g \in G$, then for any $k_1, k_2 \in O(m)$,
$$
(k_1gk_2)(k_1gk_2)' = k_1gk_2 k_2'g'k_1' = k_1gg'k_1' = k_1Xk_1' \eqdist X.
$$
Therefore, for any orthogonally invariant random matrix $X \in \P_m$, we can identify $X$ with a corresponding $g \in G$ such that $g$ is left- and right-invariant under $O(m)$.

Suppose that the random matrix $X \in \P_m$ is orthogonally invariant and $\E(X)$ is finite.  Since $X \eqdist kXk'$ for all $k \in O(m)$ then $\E(X) = \E(kXk') = k \E(X) k'$.  It follows from Schur's Lemma \cite[p.~315]{Shilov} that $\E(X) = c I_m$ for some constant $c > 0$.  

Given orthogonally invariant random matrices $X_1, X_2 \in \P_m$, the convolution product, $X_1 \circ X_2$, is defined by convolving the corresponding random matrices $g_1, g_2 \in G$ and then identifying the outcome in $\P_m$ {\it via} the natural map from $G/O(m) \to \P_m$.  Concretely, if $X_1 = g_1g_1'$ and $X_2 = g_2g_2'$ for orthogonally bi-invariant random matrices $g_1, g_2 \in G$, then $X_1 \circ X_2 = (g_1g_2) (g_1g_2)' = g_1 g_2 g_2' g_1'$.  

To express the distribution of $X_1 \circ X_2$ in terms of a product of functions of $X_1$ and $X_2$, we recall that by polar coordinates on matrix space, $g_1 = kX_1^{1/2}$ for some $k \in O(m)$.  Hence, by the orthogonal invariance of the random matrix $g$, we have $X_1 \circ X_2 = g_1g_2g_2'g_1' = g_1X_2g_1' = kX_1^{1/2} X_2 (kX_1^{1/2})' = kX_1^{1/2} X_2 X_1^{1/2}k'$.  Since $X_1$ and $X_2$ are orthogonally invariant, then replacing $X_j$ by $k'X_j k$, $j=1,2$ and using the identity $(kXk')^{1/2} = kX^{1/2}k'$ for any positive definite matrix $X$, we obtain 
$$
kX_1^{1/2} X_2 X_1^{1/2}k' \eqdist k k'X_1^{1/2}k k' X_2 k k' X_1^{1/2}kk' = X_1^{1/2} X_2 X_1^{1/2}.
$$
Consequently, 
\begin{equation}
\label{eq:composition}
X_1 \circ X_2 \eqdist X_1^{1/2} X_2 X_1^{1/2}.
\end{equation}

Although the binary operation $(A,B) \mapsto A^{1/2} B A^{1/2}$, $A, B \in \P_m$, is not associative, the following result shows that associativity holds in distribution for mutually independent, orthogonally invariant random matrices.

%\smallskip

\begin{lemma}
\label{K_invariance_assoc}
Suppose that the random matrices $X_1, X_2, X_3 \in \P_m$ are mutually independent and orthogonally invariant.  Then $(X_1 \circ X_2) \circ X_3 \eqdist X_1 \circ (X_2 \circ X_3)$.
\end{lemma}

\Pro 
By (\ref{eq:composition}),
\begin{align*}
X_1 \circ (X_2 \circ X_3) &\eqdist  X_1^{1/2} (X_2 \circ X_3) X_1^{1/2} \\
&\eqdist  X_1^{1/2} X_2^{1/2} X_3 X_2^{1/2} X_1^{1/2} = \big( X_2^{1/2} X_1^{1/2} \big)' X_3 \big( X_2^{1/2} X_1^{1/2}\big).
\end{align*}
Since $X_1^{1/2} X_2 X_1^{1/2} = (X_2^{1/2} X_1^{1/2})' (X_2^{1/2} X_1^{1/2})$ then, by polar coordinates on matrix space, there exist $k \in O(m)$ such that $X_2^{1/2} X_1^{1/2} = k (X_1^{1/2} X_2 X_1^{1/2})^{1/2}$.  Then,
\begin{align*}
X_1 \circ (X_2 \circ X_3) &\eqdist \big(  X_1^{1/2} X_2 X_1^{1/2} \big)^{1/2}  k' X_3 k \big(  X_1^{1/2} X_2 X_1^{1/2} \big)^{1/2}  \\
&\eqdist \big(X_1^{1/2} X_2 X_1^{1/2} \big)^{1/2} X_3 \big(  X_1^{1/2} X_2 X_1^{1/2} \big)^{1/2} \\
&\eqdist \big(X_1 \circ X_2 \big)^{1/2} X_3 \big(  X_1 \circ X_2 \big)^{1/2} \\
&\eqdist (X_1 \circ X_2) \circ X_3.
\end{align*}
The proof now is complete.
$\qed$

%\smallskip

\begin{lemma}
\label{K_invariance_commut}
Suppose that the random matrices $X_1, X_2 \in \P_m$ are independent and orthogonally invariant.  Then $X_1 \circ X_2 \eqdist X_2 \circ X_1$.

More generally, if the random matrices $X_1,\ldots,X_n \in \P_m$ are mutually independent and orthogonally invariant then $X_1 \circ \cdots \circ X_n \eqdist X_{\tau(1)} \circ \cdots \circ X_{\tau(n)}$ for all permutations $\tau \in \mathfrak{S}_n$, the symmetric group on $n$ symbols.
\end{lemma}

\Pro 
Since $X_1^{1/2} X_2 X_1^{1/2}$ and $X_2^{1/2} X_1 X_2^{1/2}$ have the same spectrum then there exists $k \in O(m)$ such that $X_1^{1/2} X_2 X_1^{1/2} = k X_2^{1/2} X_1 X_2^{1/2} k'$; therefore, 
$$
X_1 \circ X_2 \eqdist X_1^{1/2} X_2 X_1^{1/2} = k X_2^{1/2} X_1 X_2^{1/2} k' \eqdist k (X_2 \circ X_1)k'.
$$
Since the distributions of $X_1$ and $X_2$ are orthogonally invariant then it follows that
\begin{align*}
k(X_2 \circ X_1)k' &\eqdist k(k'X_2k \circ k'X_1k) k' \\
&= k(k'X_2k)^{1/2} k'X_1k (k'X_2k)^{1/2} k' \\
&= k k' X_2^{1/2}k k' X_1 k k'X_2^{1/2} k k' \\
&= X_2^{1/2} X_1 X_2^{1/2} \\
&\eqdist X_2 \circ X_1.
\end{align*}
Therefore, $X_1 \circ X_2 \eqdist X_2 \circ X_1$.

Finally, the proof that the distribution of $X_1 \circ \cdots \circ X_n$ is invariant under permutation of the $X_i$ is obtained by using the property that each permutation can be expressed as a product of transpositions together with the previously-proved commutativity and associativity properties.  
$\qed$

\smallskip

For $A \in \P_m$, denote by $\lambda_1(A),\ldots,\lambda_m(A)$ the eigenvalues of $A$.  The {\em Riemannian metric} on $\P_m$ is defined as 
$$
%\label{Riemannianmetric}
d_R(A,B) = \Big(\sum_{j=1}^m \big[\log \lambda_j(A^{-1/2}BA^{-1/2})\big]^2\Big)^{1/2},
$$
$A, B \in \P_m$.  The function $d_R$ is a genuine metric since, for all $A, B, C \in \P_m$: $d_R(A,B) \ge 0$, $d_R(A,B) = d_R(B,A)$, and 
\begin{equation}
\label{triangle_ineq}
d_R(A,C) \le d_R(A,B) + d_R(B,C).
\end{equation}
The quantity $d_R(A,B)$ is the distance along the shortest geodesic path, on the cone, starting at $A$ and ending at $B$ \cite{Ito,Moakher,Siegel}.  Since the matrices $A^{-1/2}BA^{-1/2}$ and $A^{-1}B$ have the same spectrum, then we also have an alternative formula, 
\begin{equation}
\label{Riemannianmetric2}
d_R(A,B) = \Big(\sum_{j=1}^m \big[\log \lambda_j(A^{-1}B)\big]^2\Big)^{1/2},
\end{equation}
from which it follows that $d_R$ is orthogonally bi-invariant: $d_R(kAk',kBk') = d_R(A,B)$ for all $k \in O(m)$ and $A, B \in \P_m$.

\begin{proposition}
\label{prop:invariancemetrics}
Let $A, B, X \in \P_m$ be orthogonally invariant random matrices.  Then 
$$
d_R(X \circ A, X \circ B) \eqdist d_R(A,B) \eqdist d_R(A \circ X, B \circ X).
$$
\end{proposition}

\Pro
By \eqref{eq:composition}, there exist $k_1, k_2 \in O(m)$ such that $X \circ A = k_1 X^{1/2} A X^{1/2}k_1'$ and $X \circ B = k_2 X^{1/2} B X^{1/2} k_2'$.  Therefore 
\begin{align*}
d_R(X \circ A, X \circ B) &= d_R(k_1 X^{1/2} A X^{1/2}k_1',k_2 X^{1/2} B X^{1/2} k_2') \\
&= d_R(k X^{1/2} A X^{1/2}k',X^{1/2} B X^{1/2}).
\end{align*}
where $k = k_2' k_1$.  By polar coordinates on matrix space, there exists $k_3 \in O(m)$ such that $k X^{1/2} = X^{1/2} k_3$; hence $X^{1/2}k' = k_3'X^{1/2}$.  Therefore, 
\begin{align*}
d_R(k X^{1/2} A X^{1/2}k',X^{1/2} B X^{1/2}) &= d_R(X^{1/2} k_3 A k_3' X^{1/2},X^{1/2} B X^{1/2}) \\
&\eqdist d_R(X^{1/2} A X^{1/2},X^{1/2} B X^{1/2}) \\
&= d_R(A,B),
\end{align*}
where the equality in distribution follows by the orthogonal invariance of $A$ and the latter equality follows from \eqref{Riemannianmetric2}.  

Again by \eqref{eq:composition} and \eqref{Riemannianmetric2}, we have 
\begin{align*}
d_R(A \circ X, B \circ X) &\eqdist d_R(A^{1/2} X A^{1/2}, B^{1/2} X B^{1/2}) \\
&= \Big(\sum_{j=1}^m \big[\log\lambda_j \big( (A^{1/2} X A^{1/2})^{-1} B^{1/2} X B^{1/2}\big)\big]^2 \Big)^{1/2}.
\end{align*}
Since $ A^{1/2} X A^{1/2} $ and $ X^{1/2} A X^{1/2} $ have the same spectrum then there exists $k_1 \in O(m)$ such that $A^{1/2} X A^{1/2} = k_1 X^{1/2} A X^{1/2} k'_1$; and similarly, there exists $k_2 \in O(m)$ such that $B^{1/2} X B^{1/2} = k_2 X^{1/2} B X^{1/2} k'_2$.  Therefore,
\begin{align*}
d_R(A \circ X, B \circ X) &= \Big(\sum_{j=1}^m \big[\log\lambda_j \big( (k_1 X^{1/2} A X^{1/2} k'_1)^{-1} k_2 X^{1/2} B X^{1/2} k'_2\big)\big]^2 \Big)^{1/2} \\
&= \Big( \sum_{j=1}^m \big[ \log\lambda_j ( k_1 X^{-1/2} A^{-1} X^{-1/2} k'_1 k_2 X^{1/2} B X^{1/2} k'_2 )\big]^2 \Big)^{1/2} \\
&= \Big(\sum_{j=1}^m \big[\log\lambda_j(k X^{-1/2} A^{-1} X^{-1/2} k' X^{1/2} B X^{1/2}) \big]^2 \Big)^{1/2}.
\end{align*}
where $k = k'_2 k_1$.  By an earlier argument, there exists $k_3 \in O(m)$ such that $kX^{-1/2} = X^{-1/2}k_3$, and hence $X^{-1/2}k' = k_3'X^{-1/2}$.  Therefore 
\begin{align*}
d_R(A \circ X, B \circ X) &= \Big(\sum_{j=1}^m \big[\log\lambda_j(X^{-1/2}k_3 A^{-1} k_3' X^{-1/2} X^{1/2} B X^{1/2}) \big]^2 \Big)^{1/2} \\
&= \Big(\sum_{j=1}^m \big[\log\lambda_j(k_3 A^{-1} k_3' B) \big]^2 \Big)^{1/2} \\
&\eqdist \Big(\sum_{j=1}^m \big[\log\lambda_j(A^{-1} B) \big]^2 \Big)^{1/2} \\
&= d_R(A,B),
\end{align*}
where the equality in distribution follows from the orthogonal invariance of $A$.  
$\qed$

\begin{remark}{\rm 
Many metrics on the cone $\P_m$ have been studied in the literature, cf., \cite{Liverani,Molnar}, and it is possible to extend the results in this section to some of those metrics.  We mention, in particular, the {\em Thompson metric}, 
\begin{align*}
d_T(A,B) &:= \log \max\{\lambda_{\max}(A^{-1/2}BA^{-1/2}),\lambda_{\max}(B^{-1/2}AB^{-1/2})\} \\
&= \log \max\{\lambda_{\max}(A^{-1}B),1/\lambda_{\min}(A^{-1}B)\},
\end{align*}
$A, B \in \P_m$, where $\lambda_{\max}(A)$ is the largest eigenvalue of $A \in \P_m$ \cite{Ito}.  It can be shown that Proposition \ref{prop:invariancemetrics} holds for $d_T$, also that 
$$
d_T(A,B) \le d_R(A,B) \le \sqrt{m} \, d_T(A,B),
$$
so that each probabilistic inequality obtained in Section \ref{sec_gen_HJ_ineqs} for $d_R$ is equivalent to a corresponding inequality for $d_T$.  
}\end{remark}

\section{Hoffmann-J{\o}rgensen inequalities on \texorpdfstring{$\boldsymbol{\P_m}$}{P(m,R)}}
\label{sec_gen_HJ_ineqs}

Let $\{X_n, n \in \N\}$ be a sequence of mutually independent, and orthogonally invariant, random matrices in $\P_m$.  Define the partial products, $S_j = X_1 \circ \cdots \circ X_j$, $j = 1,2,3,\ldots$; equivalently,
\begin{equation}
\label{Sj_matrices}
S_1 = X_1, \qquad S_{j+1} = S_j \circ X_{j+1} \equiv S_j^{1/2} X_{j+1} S_j^{1/2}, \quad j \ge 1.
\end{equation}
As we noted before, the binary relation $A \circ B$ is not commutative; however, by Lemma \ref{K_invariance_commut}, $X_1 \circ X_2 \eqdist X_2 \circ X_1$ because the random matrices $X_1$ and $X_2$ are independent and orthogonally invariant.  Therefore, we may apply the results derived by Khare and Rajaratnam \cite{Khare} (cf., \cite{Khare2}) for commutative semigroups.  

We now apply \cite[Theorem A]{Khare} with $\G$ as the collection of orthogonally invariant random matrices on $\P_m$, $\G$ being endowed with the binary operation given in \eqref{eq:composition}, and $d_\G$ being the Riemannian metric $d_R$.  For $n = 1,2,3,\ldots$ define 
\begin{equation}
\label{Mn_variables}
M_n = \max\{d_R(I_m,X_1),\ldots,d_R(I_m,X_n)\}  
\end{equation}
and
\begin{equation}
\label{Un_variables}
U_n = \max\{d_R(I_m,S_1),\ldots,d_R(I_m,S_n)\} 
\end{equation}
and let $\delta_{ij}$ denote Kronecker's delta.  We also use the notation $\chi(\mathscr{S})$ to equal $1$ if a statement $\mathscr{S}$ is valid and $0$ if it is not.  

\begin{theorem} {\rm (\cite[Theorem A]{Khare})} \, 
\label{thm_HJ1}
Let $\{X_n, n \in \N\}$ be independent, identically distributed (i.i.d.), orthogonally invariant random matrices in $\P_m$.  For $t_0, t_1, \ldots, t_l \in [0,\infty)$ and integers $l,n_1,\ldots,n_l \in \N$, let $n_{\bigcdot} = n_1+\cdots+n_l$, and define 
\begin{equation}
\label{I0_set}
\mathcal{I}_0 = \Big\{1 \le i \le l: [\Pr(U_n \le t_i)]^{n_i - \delta_{i1}} \le \frac{1}{n_i!}\Big\}.
\end{equation}
If $n_{\bigcdot} \le n+1$ then
\begin{multline}
\label{HJ1}
\Pr \Big(U_n > (n_{\bigcdot} - 1)t_0 + (2n_1-1) t_1 + 2 \sum_{j=2}^l n_j t_j\Big) \\
\le \Pr(M_n > t_0) + [\Pr(U_n \le t_1)]^{\chi(1 \notin \mathcal{I}_0)} \cdot \prod_{j \in \mathcal{I}_0} [\Pr(U_n > t_j)]^{n_j} \cdot \prod_{j \notin \mathcal{I}_0} \frac{1}{n_j!} \left(\frac{\Pr(U_n > t_j)}{\Pr(U_n \le t_j)}\right)^{n_j}.
\end{multline}

Further, let $Y_{(1)} < \cdots < Y_{(n)}$ be the order statistics of $d_R(I_m,X_1),\ldots,d_R(I_m,X_n)$.  Then the inequality \eqref{HJ1} can be strengthened by replacing $\Pr(M_n > t_0)$ by
$$
\Pr \bigg(\sum_{j=n-n_{\bigcdot}+2}^n Y_{(j)} > (n_{\bigcdot}-1)t_0 \bigg).
$$
\end{theorem}

%\medskip

\begin{remark} {\rm
Note that \eqref{HJ1} can also be written as
\begin{multline}
\label{HJ2}
\Pr \Big(U_n >  (n_{\bigcdot} - 1)t_0 + (2n_1-1) t_1 + 2 \sum_{j=2}^l n_j t_j\Big) \\
\le \Pr(M_n > t_0) + [\Pr(U_n \le t_1)]^{\chi(1 \notin \mathcal{I}_0)} \cdot \frac{\prod_{j=1}^l [\Pr(U_n > t_j)]^{n_j}}{\prod_{j \notin \mathcal{I}_0} n_j! \, [\Pr(U_n \le t_j)]^{n_j}}.
\end{multline}
}\end{remark}

%\smallskip

\begin{example} {\rm
Similar to remarks in \cite[p.~4104]{Khare}, if in Theorem \ref{thm_HJ1} we set $l = 2$, $n_1 = n_2 = 1$, and $t_1 = t_2 = t$ then by \eqref{I0_set},   
$\mathcal{I}_0 = \big\{1 \le i \le 2: [\Pr(U_n \le t)]^{1 - \delta_{i1}} \le 1\big\} = \{1,2\}$, and \eqref{HJ1} reduces to the inequality 
$$
%\label{ImIml=2}
\Pr(U_n > t_0 + 3t) \le \Pr(M_n > t_0) + [\Pr(U_n > t)]^2,
$$
which is analogous to classical Hoffmann-J{\o}rgensen inequalities \cite[p.~4101]{Khare}.
}\end{example}

%\smallskip

We prove in Appendix \ref{appendix} that if $\E(\|S_n\|_{\textrm{sp}}) < \infty$ for all $n \ge 1$ then the sequence $\{S_n, n \ge 1\}$ is a martingale and the sequences $\{M_n, n \ge 1\}$ and $\{U_n, n \ge 1\}$ are submartingales.  However, these martingale properties lead only to Markov's inequalities for $\{M_n, n \ge 1\}$ and $\{U_n, n \ge 1\}$ and in any case a martingale approach cannot lead to explicit results for specific choices of the distribution of $\{X_n, n \ge 1\}$.  Therefore we develop in Section \ref{sec_wishart} an approach that produces detailed inequalities for $\Pr(M_n > t)$ and $\Pr(U_n > t)$ in the Wishart case.  

In the following result we adopt the convention that an empty product of matrices is the identity matrix.  

%\smallskip

\begin{proposition}
\label{prop_productofXs}
Let $\{X_n, n \in \N\}$ be a sequence of i.i.d.~orthogonally invariant random matrices in $\P_m$.  Then for all $n \ge 1$, 
\begin{equation}
\label{productofXs}
S_{n+1} \eqdist X_{n+1} \circ S_n \eqdist X_{n+1}^{1/2} X_n^{1/2} \cdots X_2^{1/2} X_1 X_2^{1/2} \cdots X_n^{1/2} X_{n+1}^{1/2}.
\end{equation}
\end{proposition}

\Pro 
First, it follows from Lemma \ref{K_invariance_commut} and induction on $n$ that $S_n$ is orthogonally invariant.  

By Lemmas \ref{K_invariance_assoc} and \ref{K_invariance_commut}, the composition of independent, orthogonally invariant random matrices is commutative and associative {\it in distribution}.  Since 
\begin{align*}
S_{n+1} &= S_n \circ X_{n+1} \\
&= (S_{n-1} \circ X_n) \circ X_{n+1} \\
&= ((S_{n-2} \circ X_{n-1}) \circ X_n) \circ X_{n+1} \\
& \qquad \vdots \\
&= ((((X_1 \circ X_2) \circ X_3) \circ \cdots \circ X_{n-1}) \circ X_n)\circ X_{n+1},
\end{align*}
then it follows that 
$S_{n+1} \eqdist X_1 \circ X_2 \circ \cdots \circ X_{n+1} \eqdist X_{n+1} \circ X_n \circ \cdots \circ X_1$.  Therefore, 
$S_{n+1} \eqdist X_{n+1} \circ (X_n \circ \cdots \circ X_1) \eqdist X_{n+1} \circ S_n$.  
This proves the first equality (in distribution) in \eqref{productofXs}.  

To prove the second equality in distribution in \eqref{productofXs}, we will use induction on $n$.  Since $X_1$ and $X_2$ are i.i.d.~and $S_1 = X_1$ then 
$$
S_2 = S_1 \circ X_2 = X_1 \circ X_2 \eqdist X_2 \circ X_1 = X_2 \circ S_1 = X_2^{1/2} X_1 X_2^{1/2}.
$$
This establishes the second part of \eqref{productofXs} for $n = 1$.

For the induction step, suppose that \eqref{productofXs} holds for $n = j-1$, where $j \ge 2$.  By \eqref{Sj_matrices} and the first equality (in distribution) in \eqref{productofXs}, 
\begin{align*}
S_{j+1} &= S_{j} \circ X_{j+1} \\
&\eqdist X_{j+1} \circ S_{j} \\
&\eqdist X_{j+1} \circ (X_{j}^{1/2} X_{j-1}^{1/2} \cdots X_2^{1/2} X_1 X_2^{1/2} \cdots X_{j-1}^{1/2} X_{j}^{1/2}) \\
&\eqdist X_{j+1}^{1/2} X_{j}^{1/2} X_{j-1}^{1/2} \cdots X_2^{1/2} X_1 X_2^{1/2} \cdots X_{j-1}^{1/2} X_{j}^{1/2} X_{j+1}^{1/2}.
\end{align*}
Therefore \eqref{productofXs} holds for $n = j$.  This completes the proof by induction. 
$\qed$

\section{The random walk with Wishart matrices}
\label{sec_wishart}

In order for the inequalities \eqref{HJ1} and \eqref{HJ2} to be applied in practice, it will be necessary to obtain upper bounds for the terms $P(M_n > t)$, $P(U_n \le t)$, and $P(U_n > t)$ appearing on the right-hand sides of those inequalities.  In this section we obtain, for the case in which the random matrices $\{X_n, n \in \N\}$ are Wishart-distributed, upper bounds on $P(M_n > t)$, $P(U_n \le t)$, and $P(U_n > t)$.  

For $a > \tfrac12(m-1)$, the multivariate gamma function is defined as 
$$
%\label{Gamma_m}
\Gamma_m(a) = \pi^{m(m-1)/4} \prod_{j=1}^m \Gamma(a - \tfrac12(j-1)).
$$
For $x = (x_{ij}) \in \P_m$, let $\dd x = \prod_{1 \le i < j \le m} \dd x_{ij}$ denote Lebesgue measure on $\P_m$.  A random matrix $X \in \P_m$ is said to have the \textit{Wishart distribution}, denoted $X \sim W_m(a,I_m)$, if the probability density function of $X$ relative to Lebesgue measure on the cone is 
\begin{equation}
\label{wishart_pdf}
w(x) = \frac{1}{2^{ma}\,\Gamma_m(a)} (\det x)^{a-\frac12(m+1)} \exp(-\tfrac12 \tr x), \qquad x \in \P_m.
\end{equation}

Throughout this section, we suppose that $X_1,\ldots,X_n$ is a random sample from the Wishart distribution $W_m(a,I_m)$.  Note that once we have obtained an upper bound, $\Pr(U_n > t) \le b$ with $b < 1$, then it follows that 
$$
\frac{1}{\Pr(U_n \le t)} = \frac{1}{1 - \Pr(U_n > t)} \le \frac{1}{1-b};
$$
so we only need upper bounds for $\Pr(M_n > t)$, $\Pr(U_n > t)$, and $\Pr(U_n \le t)$.

\subsection{A bound on the distribution function of \texorpdfstring{$\boldsymbol{M_n}$}{Mn}}
\label{subsec_probs_Mn}

Define the normalizing constant, 
\begin{equation}
\label{cm_psi2}
c_m = \frac{\pi^{{m^2}/2}}{2^{ma} \, m! \, \Gamma_m(a) \, \Gamma_m(m/2)}.
\end{equation}
For $u \le v$ and $i,j=1,\ldots,m$, define 
\begin{align*}
F_{a+j}(u,v) &= \int_u^v t^{a+j-\tfrac12(m+1)} \, \exp(-t/2) \dd t, \\
F_{a+i,a+j}(u,v) &= \int_u^v t^{a+j-\tfrac12(m+1)} \, \exp(-t/2) \, F_{a+i}(u,t) \dd t,
\end{align*}
and 
$$
r_{ij}(u,v) = F_{a+i-1,a+j-1}(u,v) - F_{a+j-1,a+i-1}(u,v).
$$
Define the $m \times m$ matrix $R = \big(r_{ij}(u,v)\big)_{i,j=1,\ldots,m}$; also, for $l=1,\ldots,m$, define the $(m-1) \times (m-1)$ principal submatrices 
$$
R_l = \big(r_{ij}(u,v)\big)_{i,j=1,\ldots,l-1,l+1,\ldots,m},
$$
and the corresponding Pfaffians 
$$
H(u,v) = [\det(R)]^{1/2} \ \ \hbox{and} \ \ 
H_l(u,v) = [\det(R_l)]^{1/2}.
$$
We now have the following result.  

\begin{proposition}
\label{prop_Mn_bound}
For $t > 0$, 
\begin{equation}
\label{Mn_surviv_func_bound}
\Pr(M_n > t) \le 1 - \big[m! \, c_m \, \rho_m\big(\exp(-m^{-1/2}t),\exp(m^{-1/2}t)\big)\big]^n,
\end{equation}
where 
\begin{equation}
\label{rho_function}
\rho_m(u,v) =
\begin{cases}
H(u,v), & \hbox{if } m \hbox{ is even} \\
\operatornamewithlimits{\sum}\limits_{j=0}^{m-1} (-1)^j \, F_{a+j}(u,v) \, H_{j+1}(u,v), & \hbox{if } m \hbox{ is odd}
\end{cases}.
\end{equation}
\end{proposition}

\Pro
Since $X_1,\ldots,X_n$ are i.i.d.~then
$$
\Pr(M_n \le t) = \Pr\big(\max\big\{d_R(I_m,X_1),\ldots,d_R(I_m,X_n)\big\} \le t \big) 
= \big[\Pr\big(d_R(I_m,X_1) \le t\big)\big]^n.
$$
By \eqref{Riemannianmetric2},
$$
\Pr\big(d_R(I_m,X_1) \le t\big) = \Pr\Big(\sum_{j=1}^m [\log \lambda_j(X_1)]^2 \le t^2\Big),
$$
and since
$$
\sum_{j=1}^m [\log \lambda_j(X_1)]^2 \le m \cdot \max\{[\log \lambda_1(X_1)]^2,\ldots,[\log \lambda_m(X_1)]^2\}
$$
then it follows that 
\begin{align}
\label{Mn_bound_ineq}
\Pr\big(d_R(I_m,X_1) \le t\big) &\ge \Pr\big(m \cdot \max\big\{[\log \lambda_1(X_1)]^2,\ldots,[\log \lambda_m(X_1)]^2\big\} \le t^2\big) \nonumber \\
&= \Pr\bigg(\bigcap_{j=1}^m \{[\log \lambda_j(X_1)]^2 \le m^{-1}t^2\}\bigg) \nonumber \\
&= \Pr\bigg(\bigcap_{j=1}^m \{-m^{-1/2}t \le \log \lambda_j(X_1) \le m^{-1/2}t\}\bigg).
\end{align}
Let $\lambda_{\min}(X_1)$ and $\lambda_{\max}(X_1)$ denote the smallest and largest eigenvalues, respectively, of $X_1$.  Solving the inequalities in \eqref{Mn_bound_ineq} for $\lambda_j(X_1)$, we obtain 
\begin{align*}
\Pr(d_R(I_m,X_1) \le t) &\ge \Pr\bigg(\bigcap_{j=1}^m \{\exp(-m^{-1/2}t) \le \lambda_j(X_1) \le \exp(m^{-1/2}t)\}\bigg) \\
&= \Pr\big(\exp(-m^{-1/2}t) \le \lambda_{\min}(X_1) < \lambda_{\max}(X_1) \le \exp(m^{-1/2}t)\big),
\end{align*}
Therefore,
\begin{align}
\label{Mn_surviv_func_bound2}
\Pr(M_n > t) &= 1 - \Pr(M_n \le t) \nonumber \\
&= 1 - \big[\Pr\big(d_R(I_m,X_1) \le t\big)\big]^n \nonumber \\
&\le 1 - \big[\Pr\big(\exp(-m^{-1/2}t) \le \lambda_{\min}(X_1) < \lambda_{\max}(X_1) \le \exp(m^{-1/2}t) \big)\big]^n.
\end{align}
From a result of Krishnaiah and Chang \cite[Theorem~3.2]{Krishnaiah} for the joint distribution function of $\{\lambda_{\min}(X_1),\lambda_{\max}(X_1)\}$, we obtain the explicit, computable expression, 
$$
\Pr\big(u \le \lambda_{\min}(X_1) < \lambda_{\max}(X_1) \le v\big) = m! \, c_m \, \rho_m(u,v),
$$
where $c_m$ and $\rho_m(u,v)$ are defined in \eqref{cm_psi2} and \eqref{rho_function}, respectively.  Applying this expression to \eqref{Mn_surviv_func_bound2}, we obtain \eqref{Mn_surviv_func_bound}.  
$\qed$

\subsection{Bounds on the distribution function of \texorpdfstring{$\boldsymbol{U_n}$}{Un}}
\label{subsec_bounds_for_Un}

It is more difficult to derive bounds for the cumulative distribution function of $U_n$ than for $M_n$, the reason being that $M_n$ is defined in terms of the i.i.d.~random matrices $X_1,\ldots,X_n$ whereas $U_n$ is defined in terms of the dependent random matrices $S_1,\ldots,S_n$.  

For $0 < \theta < 1$, define 
\begin{equation}
\label{I_theta}
I_m(\theta) = (1-\theta)^{-m[(2a-m-1)(1-\theta)+m]} \prod_{j=1}^m j! \, \Gamma\big((2a-m-1)(1-\theta)+j\big).
\end{equation}
Also let 
\begin{equation}
\label{D_set}
\mathcal{D} = \{(v,\theta): 0 \le v < \tfrac12 ((2a-m-1)\theta+1), \, 0 < \theta < 1\},
\end{equation}
and define 
\begin{equation}
\label{Kvtheta_function}
K(v,\theta) = \int_0^\infty \lambda^{(2a - m-1)\theta} \exp(2v |\log \lambda| - \theta\lambda) \dd\lambda,
\end{equation}
and 
\begin{equation}
\label{Gvtheta_function}
G(v,\theta) = 2 \log c_m + \log I_m(\theta) + m \log K(v,\theta),
\end{equation}
$(v,\theta) \in \mathcal{D}$, where $c_m$ is given in \eqref{cm_psi2}.  We now establish the following result.  

\medskip

\begin{proposition}
\label{prop_Un_bound}
For $t > 0$, 
\begin{equation}
\label{Un_surviv_func_bound}
\Pr(U_n > t) \le \ n \exp\Big[\inf_{(v,\theta) \in \mathcal{D}} \big(-vt + \tfrac12 n G(v,\theta)\big)\Big].
\end{equation}
\end{proposition}

\Pro
By the exponential Markov (or Chernoff) inequality, 
\begin{equation}
\label{eq_Un_Markov}
\Pr(U_n > t) = \Pr\big(\exp(v U_n) > \exp(vt)\big) \le \exp(-vt) \E\big(\exp(v U_n)\big)
\end{equation}
for $t, v > 0$.  To obtain a bound for $\E\big(\exp(v U_n)\big)$ we apply an approach from the area of concentration inequalities, \cite[p.~375]{Boucheron}, \cite{Dasarathy}.  For $j=1,\ldots,n$, consider $\xi_j(v) = \E \exp\big(v \: d_R(I_m,S_j)\big)$, $v \ge 0$, the moment-generating function of $d_R(I_m,S_j)$.  For $j \ge 2$, we have 
\begin{align*}
\xi_j(v) &= \E \exp\big(v \: d_R(I_m,S_j)\big) \\
&= \E \exp\big(v \: d_R(I_m,X_j^{1/2}S_{j-1}X_j^{1/2})\big)
 = \E \exp\big(v \: d_R(X_j^{-1},S_{j-1})\big).
\end{align*}
By the triangle inequality \eqref{triangle_ineq}, we have
\begin{align*}
d_R(X_j^{-1},S_{j-1}) &\le d_R(X_j^{-1},I_m) + d_R(I_m,S_{j-1}) \\
&= d_R(I_m,X_j) + d_R(I_m,S_{j-1}).
\end{align*}
Exponentiating the latter inequality, taking expectations, applying the independence of $X_j$ and $S_{j-1}$, and the i.i.d. property of $X_1,\ldots,X_n$, we obtain 
\begin{align*}
\xi_j(v) &\le \E \exp\big(v \: d_R(I_m,X_j)\big) \cdot \E \exp\big(v \: d_R(I_m,S_{j-1})\big) \\
&= \E \exp\big(v \: d_R(I_m,X_1)\big) \cdot \E \exp\big(v \: d_R(I_m,S_{j-1})\big) \\
&= [\xi_1(v)] \xi_{j-1}(v). 
\end{align*}
Iterating this inequality for $j=2,\ldots,n$ we obtain 
\begin{equation}
\label{xi_j_bound}
\xi_j(v) \le [\xi_1(v)]^j
\end{equation}
for all $j \ge 1$.  Therefore 
\begin{align}
\label{Un_mgf_bound}
\E\big(\exp(v U_n)\big) &= \E\big[\exp\big(v \cdot \max\{d_R(I_m,S_1),\ldots,d_R(I_m,S_n)\}\big)\big] \nonumber \\
&= \E \big[\max\big\{\exp\big(v d_R(I_m,S_1)\big),\ldots,\exp\big(v d_R(I_m,S_n)\big)\big\}\big] \nonumber \\
&\le \sum_{j=1}^n \E \big[\exp\big(v d_R(I_m,S_j)\big)\big] \nonumber \\
&= \sum_{j=1}^n \xi_j(v).
\end{align}
Since $U_j \ge 0$ then $\xi_j(v) \ge 1$ for all $v \ge 0$; hence $\xi_j(v) \le [\xi_1(v)]^j \le [\xi_1(v)]^n$ for all $v \ge 0$ and $j=1,\ldots,n$.  By \eqref{Un_mgf_bound}, $\E\big(\exp(v U_n)\big) \le n [\xi_1(v)]^n$; also, by \eqref{eq_Un_Markov}, 
\begin{equation}
\label{E_Un_bound}
\Pr(U_n > t) \le n \exp(-vt) [\xi_1(v)]^n,
\end{equation}
a bound which is valid for all $v > 0$ such that $\xi_1(v)$ exists.  Let $\Dom(\xi_1) = \{v \ge 0: \xi_1(v) < \infty\}$; then we obtain, by minimization with respect to $v$, the inequality 
\begin{equation}
\begin{aligned}
\label{Dasarathy_bound}
\Pr(U_n > t) &\le n \inf_{v \in \Dom(\xi_1)} \exp(-vt) [\xi_1(v)]^n \\
&\equiv n \inf_{v \in \Dom(\xi_1)} \exp\big(- vt + n \log \xi_1(v)\big).
\end{aligned}
\end{equation}
We shall show later that $[0,a-\tfrac12(m-1)) \subseteq \Dom(\xi_1)$; since we have assumed that $a > \tfrac12(m-1)$ then it follows that $\Dom(\xi_1)$ contains a non-empty open set.  

Since each $X \in \P_m$ is symmetric then, by diagonalization, $X = kDk'$ where $k \in O(m)$ and $D$ is a diagonal matrix whose diagonal entries are the eigenvalues of $X$.  It follows that $f(X) = f(D)$ for any orthogonally invariant function $f: \P_m \to \C$.  Denote by $\lambda_1,\ldots,\lambda_m$ the eigenvalues of $x \in \P_m$ and denote by $\R_+$ the positive real line; then by Anderson \cite[p.~538, Theorem 13.3.1]{Anderson}, for any integrable, orthogonally invariant function $f:\P_m(\R)\to \C$, 
\begin{equation}
\label{invar_integ}
\int_{\P_m} f(x) \dd x = \frac{\pi^{m^2/2}}{m! \, \Gamma_m(m/2)} \int_{\R_+^m} f\big(\diag(\lambda_1,\ldots,\lambda_m)\big) \prod_{1 \le i<j \le m} |\lambda_i - \lambda_j| \prod_{j=1}^m \dd \lambda_j.
\end{equation}
Recalling that 
\begin{equation}
\label{xi1}
\xi_1(v) = \E d_R(I_m,X) = \E \exp\bigg(v \Big(\sum_{j=1}^m [\log \lambda_j(X)]^2\Big)^{1/2}\bigg),
\end{equation}
then, since $X \sim W_m(a,I_m)$, it follows from \eqref{wishart_pdf}, \eqref{cm_psi2}, \eqref{invar_integ}, and \eqref{xi1} that 
\begin{align}
\label{xi1_lambdas}
\xi_1(v) &= \frac{1}{2^{ma}\Gamma_m(a)} \int_{\P_m} \exp\bigg[v \Big(\sum_{j=1}^m [\log \lambda_j(x)]^2\Big)^{1/2}\bigg] (\det x)^{a-\tfrac12(m+1)} \exp(-\tfrac12 \tr x) \dd x \nonumber \\
&= c_m \int_{\R_+^m} \exp\bigg[v \Big(\sum_{j=1}^m (\log \lambda_j)^2\Big)^{1/2}\bigg] \prod_{1 \le i<j \le m} |\lambda_i - \lambda_j| \cdot \prod_{j=1}^m \lambda_j^{a-\tfrac12(m+1)} e^{-\lambda_j/2} \dd \lambda_j.
\end{align}
Also, it is straightforward to verify that 
\begin{equation}
\label{log_inequality}
\Big(\sum_{j=1}^m (\log \lambda_j)^2\Big)^{1/2} \le \sum_{j=1}^m |\log \lambda_j|,
\end{equation}
and by substituting this inequality into \eqref{xi1_lambdas} we obtain 
\begin{equation}
\label{xi1_lambdas_2}
\xi_1(v) \le c_m \int_{\R_+^m} \exp\Big[\sum_{j=1}^m (v|\log \lambda_j| - \tfrac12 \lambda_j)\Big] \prod_{1 \le i<j \le m} |\lambda_i - \lambda_j| \prod_{j=1}^m \lambda_j^{a-\tfrac12(m+1)} \dd \lambda_j.
\end{equation}

To verify that $\Dom(\xi_1) \neq \emptyset$, we now obtain {\it via} \eqref{xi1_lambdas_2} an open interval that is contained in $\Dom(\xi_1)$.  

For each $\tau \in \mathfrak{S}_m$, the group of permutations on $m$ symbols, denote by $\mathop{\rm{sgn}}(\tau)$ the signature of $\tau$.  Recall the well-known Vandermonde determinant formula, 
$$
%\label{eq_Vandermonde}
(-1)^{m(m-1)/2} \prod_{1 \le i<j \le m} (\lambda_i - \lambda_j) = \det(\lambda_i^{j-1}) = \sum_{\tau \in \mathfrak{S}_m} \mathop{\rm{sgn}}(\tau) \prod_{j=1}^m \lambda_{\tau(j)}^{j-1};
$$
consequently, for $\lambda_1,\ldots,\lambda_m > 0$, 
\begin{equation}
\label{Vandermonde_bound}
\prod_{1 \le i<j \le m} |\lambda_i - \lambda_j| \le \sum_{\tau \in \mathfrak{S}_m} \prod_{j=1}^m \lambda_{\tau(j)}^{j-1}.
\end{equation}
Substituting the bound \eqref{Vandermonde_bound} into \eqref{xi1_lambdas_2}, and using the symmetry of the resulting integrand, we find that the integral of each term on the right-hand side of \eqref{Vandermonde_bound} is the same.  Therefore $\xi_1(v) < \infty$ if 
\begin{equation}
\label{xi_1_convergence}
\begin{aligned}
\int_{\R_+^m} \exp\Big(v \sum_{j=1}^m & |\log \lambda_j|\Big) \prod_{j=1}^m \lambda_j^{a-\tfrac12(m+1)+j-1} e^{-\frac12 \lambda_j} \dd \lambda_j \\
&\equiv \prod_{j=1}^m \int_0^\infty \exp(v |\log \lambda_j| - \tfrac{1}{2}\lambda_j) \lambda_j^{(2a-m+2j-3)/2} \dd\lambda_j < \infty.
\end{aligned}
\end{equation}
The $j$th integral in this product can be written as a sum of an integral over the interval $(0,1]$ and a second integral over $(1,\infty)$, each of which is an incomplete gamma function.  Applying well-known convergence properties of the incomplete gamma function we find that the $j$th integral converges if and only if $0 \le v < (2a-m+2j-1)/2$; and from the case in which $j=1$ we deduce that $[0,a-\tfrac12(m-1)) \subseteq \Dom(\xi_1)$.  

For $0 < \theta < 1$, write  
$$
\prod_{j=1}^m \lambda_j^{a - \tfrac12(m+1)} e^{-\tfrac12 \lambda_j} = \prod_{j=1}^m \lambda_j^{\big(a - \tfrac12(m+1)\big)\theta} e^{-\tfrac12 \theta\lambda_j} \cdot \prod_{j=1}^m \lambda_j^{\big(a - \tfrac12(m+1)\big)(1-\theta)} e^{-\tfrac12 (1-\theta)\lambda_j},
$$
and substitute the right-hand side into the integrand in \eqref{xi1_lambdas}.  By the Cauchy-Schwarz inequality, 
\begin{align*}
[\xi_1(v)]^2 &\le c_m^2 \int_{\R_+^m} \exp\bigg[2v \Big(\sum_{j=1}^m (\log \lambda_j)^2\Big)^{1/2}\bigg] \prod_{j=1}^m \lambda_j^{(2a - m-1)\theta} e^{-\theta\lambda_j} \dd\lambda_j \\
& \qquad \cdot \int_{\R_+^m} \prod_{1 \le i<j \le m} |\lambda_i - \lambda_j|^2 \cdot \prod_{j=1}^m \lambda_j^{(2a - m-1)(1-\theta)} e^{-(1-\theta)\lambda_j} \dd\lambda_j.
\end{align*}
By Selberg's integral \cite[p.~406]{Andrews}, the latter integral can be evaluated explicitly: 
\begin{align}
\label{selberg1}
I_m(\theta) &:= \int_{\R_+^m} \prod_{1 \le i<j \le m} |\lambda_i - \lambda_j|^2 \cdot \prod_{j=1}^m \lambda_j^{(2a-m-1)(1-\theta)} e^{-(1-\theta)\lambda_j} \dd\lambda_j \\
&= (1-\theta)^{-m[(2a-m-1)(1-\theta)+m]} \int_{\R_+^m} \prod_{1 \le i<j \le m} |\lambda_i - \lambda_j|^2 \cdot \prod_{j=1}^m \lambda_j^{(2a-m-1)(1-\theta)} e^{-\lambda_j} \dd\lambda_j \nonumber \\
&= (1-\theta)^{-m[(2a-m-1)(1-\theta)+m]} \prod_{j=1}^m j! \, \Gamma\big((2a-m-1)(1-\theta)+j\big), \nonumber 
\end{align}
as stated in \eqref{I_theta}.  Applying the log-inequality \eqref{log_inequality}, we obtain 
\begin{align}
\label{xi1_bound}
[\xi_1(v)]^2 &\le c_m^2 I_m(\theta) \int_{\R_+^m} \exp\bigg[2v \sum_{j=1}^m |\log \lambda_j|\bigg] \prod_{j=1}^m \lambda_j^{(2a - m-1)\theta} \exp(-\theta\lambda_j) \dd\lambda_j \nonumber \\
&= c_m^2 I_m(\theta) \bigg(\int_0^\infty \exp(2v |\log \lambda| - \theta\lambda) \lambda^{(2a - m-1)\theta} \dd\lambda\bigg)^m.
\end{align}
By the convergence properties of the incomplete gamma function, we deduce that the integral in \eqref{xi1_bound} converges if and only if $(v,\theta) \in \mathcal{D}$, the set defined in \eqref{D_set}.  Therefore, by \eqref{Kvtheta_function} and \eqref{xi1_bound}, 
\begin{equation}
\label{m_1c_mIk_bound}
\xi_1(v) \le c_m [I_m(\theta)]^{1/2} \big(K(v,\theta)\big)^{m/2},
\end{equation}
$(v,\theta) \in \mathcal{D}$, and by substituting \eqref{m_1c_mIk_bound} into \eqref{Dasarathy_bound} we obtain 
$$
\Pr(U_n > t) \le n \inf_{(v,\theta) \in \mathcal{D}} \exp\big[- vt + n \log c_m + \tfrac12 n \log I_m(\theta) + \tfrac12 mn \log K(v,\theta)\big].
$$
Interchanging the infimum and the exponential function, and using the definition of $G(v,\theta)$ from \eqref{Gvtheta_function}, we obtain \eqref{Un_surviv_func_bound}.  
$\qed$

\smallskip

\begin{remark}{\rm 
(i) Proposition \ref{prop_Un_bound} implies that, for large $t$, $\Pr(U_n > t)$ converges to zero at a rate that is at least $O(e^{-ct})$, where $c > 0$.  

(ii) A bound on $\Pr(U_n > t)$ can also be obtained by using \eqref{xi_j_bound} and the resulting 3inequality, 
\begin{align*}
\E \big(\exp(v U_n)\big) \le \sum_{j=1}^n \xi_j(v) &\le \sum_{j=1}^n [\xi_1(v)]^j \\
&\le \sum_{j=1}^n [c_m [I_m(\theta)]^{1/2} \big(K(v,\theta)\big)^{m/2}]^j \\
&= \frac{\big(c_m [I_m(\theta)]^{1/2} (K(v,\theta))^{m/2}\big)^{n+1} -  1}{c_m [I_m(\theta)]^{1/2} (K(v,\theta))^{m/2} - 1} - 1.
\end{align*}
This bound is sharper than the one given in \eqref{Dasarathy_bound}, so it can be expected to lead to sharper inequalities for $\Pr(U_n > t)$.  However, the resulting numerical calculations will be more involved.  
}\end{remark}

We now consider the issue of calculating the infimum in \eqref{Un_surviv_func_bound} for given values of $m$, $n$, and $a$.  We will see that the objective function, $(v,\theta) \mapsto -vt + \tfrac12 m G(v,\theta)$, $(v,\theta) \in \mathcal{D}$, is particularly well-behaved as it is strictly convex in $\mathcal{D}$; consequently, the infimum occurs at a unique point in $\mathcal{D}$.  

\begin{proposition}
\label{prop_convexity}
(i) For $(v,\theta) \in \mathcal{D}$, the function $K(v,\theta)$, is strictly log-convex in $(v,\theta)$ jointly, strictly log-convex in $v$ for fixed $\theta$, and strictly log-convex in $\theta$ for fixed $v$.  

(ii) The function $I_m(\theta)$, $\theta \in (0,1)$, is strictly log-convex.

(iii) For $t \ge 0$, the objective function $(v,\theta) \mapsto -vt + \tfrac12 m G(v,\theta)$, $(v,\theta) \in \mathcal{D}$, is strictly convex.
\end{proposition}

\Pro
These properties follow readily from the fact that $K(v,\theta)$ and $I_m(\theta)$ can be viewed as modified Laplace transforms of positive measures.  For instance, let us write $K(v,\theta)$ in the form 
$$
K(v,\theta) = \int_0^\infty \exp\big(2v \phi_1(\lambda) - \theta\lambda + \phi_2(\lambda)\big) \dd \lambda, 
$$
where $\phi_1(\lambda) = |\log \lambda|$ and $\phi_2(\lambda) = (2a - m-1)\theta \log \lambda$.  By proceeding as in the theory of exponential families of distributions (Barndorff-Nielsen \cite[p.~103 ff.]{BN}), the standard method of applying H\"older's inequality leads to the conclusion that $K(v,\theta)$ is log-convex in $(v,\theta)$ jointly.  Moreover, the strict log-convexity of $K(v,\theta)$ follows from the fact that the function 
$$
\bigg|\begin{matrix}
\phi_1(\lambda_1) & \phi_1(\lambda_2) \\
\phi_2(\lambda_1) & \phi_2(\lambda_2) \\
\end{matrix}\bigg|^2 ,
$$
$(\lambda_1,\lambda_2) \in \R_+^2$, is strictly positive on an open subset of $\R_+^2$.  

The remaining statements in (i) can be proved similarly, (ii) can be proved analogously using the integral representation given in \eqref{selberg1}, and (iii) follows from (i) and (ii).
$\qed$

\medskip

Once the strict convexity property of the objective function has been established, it is numerically straightforward to compute in \eqref{Un_surviv_func_bound} both the value of the infimum and its point of location.  We remark that it is especially simple to compute numerical examples using, e.g., the \textsl{SageMath} software system \cite{Zimmerman}.  

Next, we obtain an upper bound for the distribution function, $\Pr(U_n \le t)$.  Denoting by $\chi^2_\alpha$ a random variable having the chi-squared distribution with degrees-of-freedom $\alpha$, we have the following result.  

\begin{proposition}
\label{prop_Un_cdf_bound}
Let $R_1,\ldots,R_n$ be i.i.d., $\chi^2_{(2a-m-1)m}$-distributed random variables.  Then, 
\begin{equation}
\label{Un_cdf_bound}
\Pr(U_n \le t) \le \min_{1 \le j \le n} \Pr\Big(\prod_{i=1}^j R_i \le m^j \exp(m^{-1/2}t)\Big),
\end{equation}
$t > 0$.  In particular, 
\begin{equation}
\label{Un_cdf_bound_j1}
\Pr(U_n \le t) \le \Pr\big(\chi^2_{(2a-m-1)m} \le m \exp(m^{-1/2}t)\big).
\end{equation}
\end{proposition}

\Pro
By Bonferroni's inequality, 
\begin{align}
\label{Un_Bonferroni}
\Pr(U_n \le t) &= \Pr\big(\max\{d_R(I,S_1),\ldots,d_R(I,S_n)\} \le t\big) \nonumber \\
&= \Pr\big(d_R(I,S_1) \le t,\ldots,d_R(I,S_n) \le t\big) \nonumber \\
&\le \min_{1 \le j \le n} \Pr\big(d_R(I,S_j) \le t\big).
\end{align}
For $j=1,\ldots,n$, 
\begin{align*}
\Pr\big(d_R(I,S_j) \le t\big) &= \Pr\big([d_R(I,S_j)]^2 \le t^2\big) 
= \Pr\Big(\sum_{i=1}^m [\log\lambda_i(S_j)]^2 \le t^2\Big).
\end{align*}
By the Cauchy-Schwarz inequality, 
\begin{align*}
\sum_{i=1}^m [\log\lambda_i(S_j)]^2 &\ge \frac{1}{m} \Big(\sum_{i=1}^m \log\lambda_i(S_j)\Big)^2 \\
&= \frac{1}{m} \Big(\log \prod_{i=1}^m \lambda_i(S_j)\Big)^2 
= \frac{1}{m} \big[\log (\det S_j)\big]^2.
\end{align*}
Therefore, 
\begin{align*}
\Pr\big(d_R(I,S_j) \le t\big) &\le \Pr\Big(\frac{1}{m} \big[\log (\det S_j)\big]^2 \le t^2\Big) 
= \Pr\big((\det S_j) \le \exp(m^{1/2}t)\big).
\end{align*}
By \eqref{productofXs}, $\det S_j \eqdist \prod_{i=1}^j (\det X_i)$, and we proceed to obtain a lower bound on the distribution of this product.  

For random variables $Q_1$ and $Q_2$, recall that $Q_1$ is said to be {\it stochastically greater than} $Q_2$, written $Q_1 \gedist Q_2$, if $\Pr(Q_1 \le t) \le \Pr(Q_2 \le t)$ for all $t \in \R$.  It is well-known that if $Q_{ij}$, $i,j=1,2$ are positive random variables such that $Q_{11} \gedist Q_{21}$ and $Q_{12} \gedist Q_{22}$ then $Q_{11} Q_{12} \gedist Q_{21} Q_{22}$.  

Since $X_j \sim W_m(a,I_m)$ then, by \cite[Theorem 4]{Gordon}, 
$$
m (\det X_j)^{1/m} \gedist \chi^2_{(2a-m-1)m}.
$$
Therefore, since $R_1,\ldots,R_n$ are i.i.d.~and $\chi^2_{(2a-m-1)m}$-distributed then it follows that 
$$
m^j (\det S_j)^{1/m} \eqdist \prod_{i=1}^j m (\det X_i)^{1/m} \gedist \prod_{i=1}^j R_i,
$$
and therefore 
\begin{align*}
\Pr\big(d_R(I,S_j) \le t\big) &\le \Pr\big(\det S_j \le \exp(m^{1/2}t)\big) \\
&= \Pr\big(m^j (\det S_j)^{1/m} \le m^j \exp(m^{-1/2}t)\big) \\
&\le \Pr\Big(\prod_{i=1}^j R_i \le m^j \exp(m^{-1/2}t)\Big).
\end{align*}
Applying \eqref{Un_Bonferroni}, we now obtain \eqref{Un_cdf_bound}.  Also, by choosing $j=1$ we obtain \eqref{Un_cdf_bound_j1}.
$\qed$

\medskip

We remark that much is known about the exact distribution of a product of independent chi-square random variables; cf.~\cite{Anderson,Gordon}.  Such results can be used to obtain sharper, but more complicated, upper bounds on the distribution function of $U_n$.  

\bigskip
\medskip

\noindent
\textbf{Acknowledgments}:  We are very grateful to Apoorva Khare for providing us with comments that led to improvements in an earlier draft of this manuscript.  We also thank an anonymous referee for reading our manuscript with great thoroughness.

\medskip

\noindent
\textbf{Declaration}:  No financial support was received to underwrite the preparation of this manuscript, and the authors have no relevant financial or non-financial interests to disclose.

\medskip

\noindent
\textbf{Data Availability Statement}:  Data sharing is not applicable to this article as no datasets were generated or analysed during the current study.

\appendix

\section{Martingale properties of \texorpdfstring{$\boldsymbol{S_n}$}{Sn}, \texorpdfstring{$\boldsymbol{M_n}$}{Mn}, and \texorpdfstring{$\boldsymbol{U_n}$}{Un}}
\label{appendix}

Let $\{Y_n, n \in \N\}$ be a sequence of random entities taking values in a finite-dimensional Euclidean space $(\R^d,\|\cdot\|)$, where $\E(\|Y_n\|) < \infty$.  The sequence $\{Y_n, n \in \N \}$ is called a {\it martingale} if $Y_n = \E(Y_{n+1} | Y_1,\ldots,Y_n)$ for all $n \ge 1$.  If each $Y_n$ is scalar-valued then the sequence $\{Y_n, n \in \N \}$ is called a {\it submartingale} if $Y_n \le \E(Y_{n+1} | Y_1,\ldots,Y_n)$ for all $n \ge 1$.  

\begin{proposition}
\label{martingale_Sn}
Let $\{X_n, n \in \N\}$ be a sequence of orthogonally invariant, i.i.d.~random matrices in $\P_m(\R)$ such that $\E(X_n) = I_n$ for all $n$.  Then the sequence $\{S_n, n \in \N\}$ in \eqref{Sj_matrices} is a martingale, and the sequences $\{M_n, n \in \N\}$ in \eqref{Mn_variables} and $\{U_n, n \in \N\}$ in \eqref{Un_variables} are submartingales.
\end{proposition}

\Pro 
Let $\S_m(\R)$ denote the vector space of $m \times m$ symmetric matrices, and endow $\S_m(\R)$ with the spectral norm, $\|X\|_{\textrm{sp}} := \max\{|\lambda_1(X)|,\ldots,|\lambda_m(X)|\}$, the maximum absolute value of the eigenvalues of $X$.  For $X \in \P_m(\R)$, it is evident that $\|X^{1/2}\|_{\textrm{sp}} = \|X\|_{\textrm{sp}}^{1/2}$.  Also, by \cite[p.~343, Theorem 5.6.2]{Horn}, this norm is submultiplicative: $\|XY\|_{\textrm{sp}} \le \|X\|_{\textrm{sp}} \cdot \|Y\|_{\textrm{sp}}$.  Therefore, 
\begin{align*}
\|S_n\|_{\textrm{sp}} \eqdist \|X_n \circ S_{n-1}\|_{\textrm{sp}} 
&= \|X_n^{1/2} S_{n-1} X_n^{1/2}\|_{\textrm{sp}} \\
&\le (\|X_n\|_{\textrm{sp}}^{1/2})^2 \cdot \|S_{n-1}\|_{\textrm{sp}} = \|X_n\|_{\textrm{sp}} \|S_{n-1}\|_{\textrm{sp}}.
\end{align*}
By iterating this inequality, or by applying \eqref{productofXs}, we obtain $\|S_n\|_{\textrm{sp}} \ledist \|X_1\|_{\textrm{sp}} \cdots \|X_n\|_{\textrm{sp}}$.  

Since $X_1,\ldots,X_n$ are i.i.d.~then it follows that 
$$
\E(\|S_n\|_{\textrm{sp}}) \le \E(\|X_1\|_{\textrm{sp}} \cdots \|X_n\|_{\textrm{sp}}) = (E \|X_1\|_{\textrm{sp}})^n.
$$
Also, since $\|X_1\|_{\textrm{sp}} < \tr X_1$ then 
$$
\E(\|X_1\|_{\textrm{sp}}) \le \E(\tr X_1) = \tr \E(X_1) = m < \infty;
$$
therefore, $\E(\|S_n\|_{\textrm{sp}}) \le m^n < \infty$.  

By \eqref{productofXs} and the independence of $X_{n+1}$ and $\{X_1,\ldots,X_n\}$, 
\begin{align*}
\E(S_{n+1} | S_1,\ldots,S_n) &= \E(S_{n}^{1/2} X_{n+1} S_n^{1/2} | S_1, \ldots, S_n) \\
&= S_{n}^{1/2} \E(X_{n+1} | S_1,\ldots,S_n) S_{n}^{1/2}
= S_{n}^{1/2} \E(X_{n+1}) S_{n}^{1/2} = S_n.
\end{align*}
Therefore, the sequence $\{S_n, n \in \N\}$ is a martingale. 

We now show that $\{M_n, n \in \N\}$ is a submartingale.  First,  
$$
\E(M_n) = \E \max_{1 \le j \le n} d_R(I_m,X_j) \le \sum_{j=1}^n \E d_R(I_m,X_j) = n \E d_R(I_m,X_1).
$$
By \eqref{log_inequality}, 
$$
\E d_R(I_m,X_1) = \E \Big(\sum_{j=1}^n [\log \lambda_j(X_1)]^2\Big)^{1/2} \le \E \sum_{j=1}^n |\log \lambda_j(X_1)|.
$$
Applying the same analysis as at \eqref{xi_1_convergence}, we find that the latter expected value is finite.  See also the functions $\xi_1(v)$ and $k(v,\theta)$ in \eqref{xi1_lambdas} and \eqref{Kvtheta_function}.  Therefore $\E(M_n) < \infty$.  

Also, it is obvious that $M_{n+1} \ge M_n$ for all $n \ge 1$.  Therefore 
$$
\E(M_{n+1} | M_1,\ldots,M_n) \ge \E(M_n | M_1,\ldots,M_n) = M_n,
$$
so $\{M_n, n \in \N\}$ is a submartingale.  

As for $\{U_n, n \in \N\}$, it follows from \eqref{E_Un_bound} that 
$$
\E(U_n) = \int_0^\infty \Pr(U_n > t) \dd t \le n v^{-1} [\xi_1(v)]^n,
$$
and therefore 
$$
\E(U_n) \le n \inf_{v \in \Dom(\xi_1)} v^{-1} [\xi_1(v)]^n < \infty.
$$
Further, $U_{n+1} \ge U_n$ for all $n \ge 1$ and therefore 
$$
\E(U_{n+1} | U_1,\ldots,U_n) \ge \E(U_n | U_1,\ldots,U_n) = U_n,
$$
so $\{U_n, n \in \N\}$ also is a submartingale.  
$\qed$

\bigskip 

By Kolmogorov's inequalities for submartingales \cite[p.~314]{Ross}, 
$$
\Pr(M_n \ge t) \le t^{-1} \E(M_n)
$$
and 
$$
\Pr(U_n \ge t) \le t^{-1} \E(U_n), 
$$
$t > 0$.  Since the sequences $\{M_n, n \in \N\}$ and $\{U_n, n \in \N\}$ are nonnegative and increasing then it follows that Kolmogorov's inequalities for $M_n$ and $U_n$ are the same as Markov's inequalities.  

\end{document}